\documentclass[11pt, a4paper]{amsart}
\newcommand{\Frc}{\mbox{frac}}
\usepackage{bbm,latexsym,amsmath,amssymb, verbatim}

\usepackage[english]{babel}
\usepackage{epsfig}
\usepackage[sectionbib]{natbib}
\usepackage{placeins}
\usepackage{float}
\newcommand\4{\hskip 1pt}

\newcommand{\beq}{\begin{equation}}
\newcommand{\beqa}{\begin{eqnarray}}
\newcommand{\eeqa}{\end{eqnarray}}
\newcommand{\eeq}{\end{equation}}
\newcommand{\bdi}{\begin{displaymath}}
\newcommand{\edi}{\end{displaymath}}
\renewcommand\Re{\operatorname{Re}}
\renewcommand\Im{\operatorname{Im}}
\newcommand{\q}{\quad}

\newcommand{\Ee}{\mathbf E}
\newcommand{\Ze}{\hbox{\rm Z\negthinspace \negthinspace Z }}
\numberwithin{equation}{section}
\renewcommand\Re{\operatorname{Re}}
\renewcommand\Im{\operatorname{Im}}

\begin{document}

\title[Student Statistic]{
An Expansion for a Discrete Non-Lattice Distribution}
\author[F. G\"otze]{Friedrich G\"otze$^1$}
\address{Faculty of Mathematics, University of Bielefeld, Germany}
\email{goetze@math.uni-bielefeld.de}

\author[W.R. van Zwet]{Willem R. van Zwet $^1$}
\address{Department of Mathematics, Leiden University, Netherlands}
\email{vanzwet@math.leidenuniv.nl}
\thanks{$^1$Research supported by CRC 701}
 
\date{September 9, 2005}

\begin{abstract}
Much is known about asymptotic expansions for asymptotically normal distributions if these distributions are either absolutely continuous or pure lattice distributions. In this paper we begin an investigation of the discrete but non-lattice case. We tackle one of the simplest examples imaginable and find that curious phenomena occur. Clearly more work is needed.
\end{abstract}

\maketitle



\section{Introduction}     

There is a voluminous literature on second order analysis of distribution functions $F_N(z) = P(Z_N\leq z)$ of statistics $Z_N = \zeta_N(X_1,X_2,\dots, X_N)$ that are functions of i.i.d. random variables $X_1,X_2,...$. The results obtained are generally refinements of the central limit theorem. Suppose  that $Z_N$ is asymptotically standard normal, that is $\sup_z |F_N(z)-\Phi (z)|\rightarrow0 $ as $N\rightarrow \infty$, where $\Phi $  denotes the standard normal distribution function. Then second order results are concerned with the speed of  this convergence, or with attempts to increase this speed by replacing the limit $\Phi $ by a series expansion $\Psi_N$ that provides a better approximation. Results of the first kind are called theorems of Berry-Esseen type and assert that for all $N$,
$$\sup_z|F_N(z)-\Phi(z)|\leq CN^{-\frac12},$$
where $C$ is a constant that depends on the particular statistic $Z_N$ and the
distribution of the variables $X_i$, but not on $N$. Such results are often
valid under mild restrictions such as the existence of a third moment of
$X_i$. The original Berry-Esseen theorem dealt with the case where $Z_N$ is a
sum of i.i.d. random variables, 
\citet*{Esseen:1942},
\citet*{Berry:1941}. 
For a more general version see 
\citet*{vanZwet:1984}.

Results of the second kind concern so-called Edgeworth expansions. These are series expansions such as 
\begin{equation} 
        \begin{split}
         &\Psi_{ N,1}(z) = \Phi(z)+\varphi(z)N^{-\frac12} Q_1(z), \q \text{or} \\ 
        &\Psi_{N,2}(z)  =\Phi(z)+\varphi(z)\left[ N^{-\frac12}  Q_1(z)+ N^{-1}Q_2(z)\right],
        \end{split}
\end{equation}  

\noindent
where $ \varphi $ is the standard normal density and $Q_1$ and $Q_2$ are polynomials depending on low moments of $X_i$ . One then shows that 
\begin{equation} 
        \begin{split}
                &\sup_z|F_N(z)-\Psi_{N,1}(z)|\leq CN^{-1}, \q \text{or}  \\  
                &\sup_z|F_N(z)-\Psi_{N,2}(z)|\leq CN^{-\frac32} .
        \end{split}
\end{equation} 

\noindent
For this type of result the restrictions are more severe. Apart from moment
assumptions one typically assumes that $Z_N$ is not a lattice random
variable. For the case where $Z_N$ is a sum of i.i.d. random variables a good
reference is 
\citep[chap. XVI]{Feller:1965/2}. 
There are numerous papers devoted to special types of statistics. For a
somewhat more general result we refer to 
\citet*{Bickel-Goetze-vanZwet:1986}
 and 
\citet*{bentkus-goetze-vanzwet:1997} .

For the case where $Z_N$ assumes its values on a lattice, say the integers, an alternative approach it to generalize the local central limit theorem and provide an expansion for the point probabilities  $P(Z_N=z)$ for values of $z$ belonging to the lattice. A typical case is the binomial distribution for which local expansions are well known. It is obvious that for the binomial distribution one can not obtain Edgeworth expansions as given in $(1.1)$  for which $(1.2)$ holds. The reason is that out of the $N$ possible values for a binomial $(N,p)$ random variable, only $cN^{\frac12}$ values around the mean $Np$ really count and each of these has probability of order $N^{-\frac12}$. Hence the distribution function has jumps of order $N^{-\frac12}$ and can therefore not be approximated by a continuous function such as given in $(1.1)$ with an error of smaller order than $N^{-\frac12}$. 

In a sense the binomial example is an extreme case where the ease of
the approach through local expansions for $P(Z_N=z)$ precludes the one through
expansions of Edgeworth type for $P(Z_N\leq z)$. In 
\citet*{Albers-Bickel-vanZwet:1976}
 the authors found somewhat to their surprise that for the Wilcoxon statistic which is a pure lattice statistic, an Edgeworth expansion with remainder $O(N^{-\frac32})$ for the distribution function is perfectly possible. In this case the statistic ranges over $N^2$ possible integer values, of which the central $N^{\frac32}$ values have probabilities of order $N^{-\frac32}$ so that one can approximate a distribution function with such jumps by a continuous one with error $O(N^{-\frac32})$.
 
On the basis of these examples one might guess that the existence of an Edgeworth  expansion with error $O(N^{-p})$ for the distribution function $F_N(z) = P(Z_N\leq z)$ would merely depend on the existence of some moments of $Z_N$ combined with the requirement that $F_N$ does not exhibit jumps of large order than $N^{-p}$. But one can envisage a more subtle problem if $F_N$ would assign a probability of larger order than $N^{-p}$ to an interval of length $N^{-p}$. Since Edgeworth expansions have bounded derivative, this would also preclude the existence of such an expansion with error $O(N^{-p})$.
  
Little seems to be known about the case where $Z_N$ has a discrete but non-lattice distribution. Examples abound if one considers a lattice random variable with expectation $0$ and standardized by dividing by its sample standard deviation. As a simple example, one could for instance consider Student's $t$-statistic $\tau_N = N^{-1/2}\sum_i X_i / \sqrt{\sum_i\left( X_i -m \right)^2/(N-1) }$  with $m= \sum_i X_i/N $ and $X_1,X_2,\dots$ i..i.d. random variables with a lattice distribution. Since we are not interested in any particular statistic, but merely in exploring what goes on in a case like this, we shall simplify even further by deleting the sample mean m and considering the statistic\smallskip

\noindent 
\begin{equation}   
W_N = \sum_{i=1}^N \frac{X_i}{\sqrt{\sum_{i=1}^N X_i^2}} ,
\end{equation}

\noindent
with \smallskip

\noindent 
\begin{equation}
 X_1,X_2,\dots i.i.d.\mbox{ with } P(X_i=-1)= P(X_i=0)= P(X_i=1)= \frac13 .
\end{equation}

We should perhaps point out that for $w>0$
\begin{align*}
P(0<\tau_N\leq w)&=P\left(0<W_N\leq \frac{\sqrt{N/(N-1)}x}{1+x^2/(N-1)}\right) \\
&= P(0<W_N\leq x + O(x(1+x^2)/N)),
\end{align*}
and since both $\tau_N$ and $W_N$ have distributions that are symmetric about the origin, Theorem 1.1 ensures that under this model, the expansions for the distributions of $\tau_N$ and $W_N$ are identical up to order $O(N^{-1})$. Hence the fact that we discuss an expansion for $W_N$ rather than for Student's statistic $\tau_N$ is not an essential difference, but merely a matter of convenience.

Notice that in (1.3) $ \sum X_i^2$ equals the number of non-zero $X_i$ and that each of these equals $-1$ or $+1$ with probability $\frac12$. Hence we may also describe the model given by (1.3) and (1.4) as follows.

For $N=1,2,\dots, T_N$ has a binomial distribution with parameters $N$ and $\frac23$. Given $T_N, S_N$ has a binomial distribution with parameters $T_N$ and $\frac12$ . Define $D_N=2 S_N-T_N$ and notice that $D_N$ and $T_N$ are either both even or both odd. We consider the statistic \smallskip

\noindent 
 \begin{equation}W_N = \begin{cases} 
                       0  \q \text{ if } T_N=0 \\
                       \frac{D_N}{\sqrt{T_N} } \q \text{ otherwise.}                                           
                    \end{cases}
        \end{equation}

\noindent
Notice that unconditionally $S_N$ has a binomial distribution with parameters $N$ and $\frac13$.

  Let $F_N$ be the distribution function of  $W_N$.  Obviously $W_N$ has
  expected value $EW_N=0$ and variance $\sigma^2 (W_N)=1$, and is
  asymptotically normal. By an appropriate Berry-Esseen theorem $\sup_x
  |F_N(x)-\Phi(x)| = O(N^{-\frac12})$, where $\Phi$ denotes the standard
  normal distribution function as before, see
  e.g. 
\citet*{Bentkus-Goetze:1996}.
 A curious thing about the distribution of $W_N$ is that $ P(W_N =0)=P(D_N=0)=P(S_N=\frac12T_N)$ which is obviously of exact order $N^{-\frac12}$, but all other point probabilities $P(W_N =w)$ for $w\not= 0$ are clearly of smaller order, and we shall see that these are actually $O(N^{-1})$. Hence the following question arises: if we remove the point probability at the origin, to what order of magnitude can we approximate the distribution function of $W_N$  by an expansion?

Let $\Frc(x)$ denote the fractional part of a number $x\geq 0$, i.e. $\Frc (x)=x-\lfloor x\rfloor $ if $\lfloor x\rfloor $ is the largest integer smaller than or equal to $x$. For $N=1,2,\dots$, define a function $\Psi_N$ on $(-\infty,\infty)$ as follows.\\ For $w\geq 0$,\smallskip

\noindent 
\begin{equation*}
\Psi_N(w) = \Phi(w) + N^{-1/2} \Lambda_N(w),\ {\rm with}
\end{equation*}          
\begin{equation}
\Lambda_N(w) = -\sqrt{ \frac{3}{2}}\, \varphi(w)  \sum_{0\leq n
           \leq N} \frac{3}{\sqrt{\pi N}} \,e^{    -\frac{9}{N}\left(n-\frac
             N3\right)^2  }  
\left(\Frc \left(  w \sqrt{2n}  \right)-\frac12 \right)
\end{equation}
and
\begin{equation*}
\Psi_N(-w)= 1-\Psi_N(w-).
\end{equation*} 
Here the argument $w-$ denotes a limit from the left at $w$. $\Psi_N$ is of bounded variation and for sufficiently large $N$ it is a probability distribution function. It has upward jump discontinuities of order $O(N^{-1})$ at points $w$ where $w\sqrt{2n}$ assumes an integer value $k\not=0$. At the origin it has a jump discontinuity of magnitude 
\begin{equation} 
        \Psi_N(0)- \Psi_N(0-) \sim \sqrt{\frac{3}{2N}}\varphi(0)=\sqrt{\frac{3}{4\pi N}}.
\end{equation}
 
\noindent
{\bf Theorem 1.1.} \q {\it As} $N\rightarrow \infty ,$
\begin{equation*}
 P(W_N=0)\sim \,\,\sqrt{\frac{3}{(4\pi N)} }
=O\left( \frac{1}{\sqrt{N}}\right),
\end{equation*}
\begin{equation}
P(W_N = w) =  \,\,O\left(\frac 1N\right),
\end{equation}
{\it for $w \not= 0$, and}          
\begin{equation*} 
\sup_{w} |F_N(w)-\Psi_N(w)| = \,\,O\left(\frac 1N\right).
\end{equation*}
{\it Moreover the distribution $\Psi_N$ - and hence $F_N$ -  assigns probability $O(N^{-1})$ to any closed interval of length $O(N^{-1})$ that does not contain the origin. }

The reader may want to compare this result with the expansion obtained in 
\citet*{Brown-Cai-DasGupta:2002} 
for the distribution of a normalized binomial variable $(Y_n-np)/\sqrt{np(1-p)}$, where $Y_n$ has a binomial distribution with parameters $n$ and $p$. For $p=1/2$ this distribution coincides with the conditional distribution of $W_N$ given $N=n$, but of course this is a pure lattice distribution.

Since $\sum_n 3(\pi N)^{-1/2} \exp\{-(9/N)(n-N/3)^2\}$ is the sum of the density of a normal $N(N/3,N/18)$ distribution taken at the integer values $n \in (-\infty,\infty)$, it is asymptotic to 1 as $N \to \infty$, and hence bounded for all $N$. Because $|\Frc(x)-1/2|\leq 1/2$ for all $x>0$, $\Lambda_N(w)$ is bounded and $N^{-1/2}\Lambda_N(w) = O(N^{-1/2})$ uniformly in $w$. At first sight there is a striking similarity between the expansion $\Psi_N$ in Theorem 1.1 and the two term Edgeworth expansion $\Psi_{N,1}$ in (1.1). However, the term $\phi(z)N^{-1/2}Q_1(z)$ of order $O(N^{-1/2})$ in the Edgeworth expansion is a skewness correction that vanishes for a symmetric distribution $F_N$. As we are dealing with a symmetric case, such a term is not present and for the continuous case the Edgeworth expansion with remainder $O(N^{-1})$ is simply $\Phi(z)$. The origin of the term $N^{-1/2}\Lambda_N(w)$ is quite different. It arises from the fact that we are approximating a discrete distribution function by a continuous one, and as such it is akin to the classical continuity correction.

To make sure that the term $N^{-1/2}\Lambda_N(w)$ is not of smaller order than $N^{-1/2}$, we shall bound $|\Lambda_N(w)|$ from below  by the absolute value 
of the following series. Assume that $N$ is divisible by $3$ and let
  \begin{equation}
        \begin{split}
\lambda_N(w):=& \sqrt{\frac 3 2}\4 \varphi(w)\sum_{k=1}^{M} \frac 1 {\pi \4 k} 
f_{N,k}\exp\bigl(- \frac{\pi^2}6\4 k^2\4  w^2\4)\4 
\sin\bigl(2 \4 \pi\4 \4 k \4 w\4 \sqrt{\frac {2 N} 3}\bigr)\\
& \, \q + O\bigl(N^{-1/2}(\log N)^5\bigr), \q M := \lfloor \log N \rfloor ,
 \end{split}
        \end{equation}
where  $f_{N,k}=1 + O((k/M)^2)$ is defined in $(3.9)$. Thus $\lambda_N(w)$ 
is a rapidly converging Fourier series, (illustrated in
Figure  1. below)  the modulus of which is larger than a positive  constant $c(w)>0$,
provided that $4\4 w\4 \sqrt{\frac{2 N} 3} $ is an odd integer.

\begin{figure}[H]          
\psfig{file=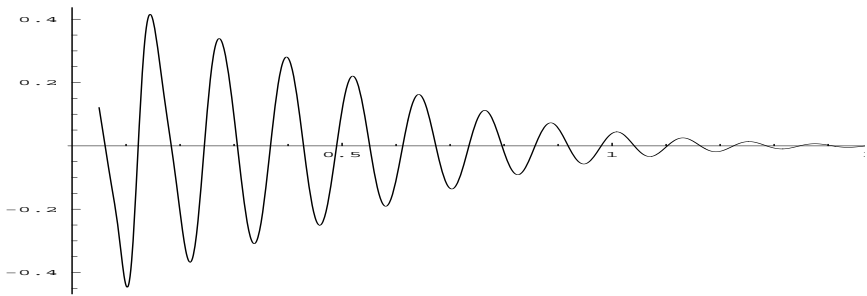,width=18cm, height=4cm}
\caption{$\lambda_{100}(w)$:\,\, $0.05 \le w \le 2.34,\,\, M=10$, $f_{100,k}:=\exp[-(k/M)^{2/3}]$\label{fig1}} 
\end{figure}  
 Hence, we shall prove
\noindent
{\bf Theorem 1.2.} {\it For any $N$ divisible by $3$, we have} 
\begin{equation}
                 \sup_{w>0} |F_N(w)-\Phi(w)| \ge  \,\,
\sup_{w\ge 1}N^{-\frac 12}|\lambda_N(w)| + O\bigl(N^{-1} (\log N)^5\bigr)
>\,\, \frac c {\sqrt N},
\end{equation}
{\it for some absolute constant $c>0$.} 

\bigskip\noindent
The proof of Theorem 1.1 is given in Section 2. In Section 3 we investigate 
the oscillatory part of $\Psi_N$ in (1.6), relating it to the
Fourier series $\lambda_N(w)$ above and thus proving  Theorem 1.2.

\noindent
{\bf Acknowledgment.}
 The authors would like to thank G.Chistyakov  for a careful reading of the
 manuscript
 and Lutz Mattner for his comments on the  current ArXiv version.
 

\section{Proof of Theorem 1.1}

The event $W_N = 0$ occurs iff $D_N = 0$.  Let $Z_1,Z_2,\dots, Z_N$ be i.i.d. random variables assuming the values $0$, $-1$ and $+1$, each with probability $\frac 13$. Then $D_N$ is distributed as  $ \sum Z_i$, which has mean $0$ and variance $\frac{2N}{3}$. By the local central limit theorem  $P(\sum Z_i =0) \sim (2\pi)^{-\frac12} \left(\frac{2N}{3} \right)^{-\frac12} = \sqrt{\frac{3}{4\pi N}}$ which proves the first statement of Theorem 1.1. Because the distribution of $W_N$ is symmetric about the origin, this implies that in the remainder of the proof we only need to consider positive values of $W_N$. Hence we suppose that $w>0$ throughout and this implies that we need only be concerned with positive values of  $D_N$ also.

Hoeffding's inequality ensures that for all $N\geq 2$,
 \[P\left(|D_N|\geq  \sqrt{6N \log N} \right)\leq   \frac{2}{N^{3}}   \]
and 
\[P\left(|T_N -2N/3|\geq \sqrt{2 N \log N}\right)  \leq  \frac{2}{N^{3}}.\] 
\noindent
Since the joint distribution of $T_N$ and $D_N$ assigns positive probability to at most $N^2$ points and events with probability $O(N^{-1})$ are irrelevant for the remainder of the proof,  we may at any point restrict attention to values $D_N=d$ and $T_N=t$ with $|d| \leq t$ and satisfying \smallskip

\noindent
          \begin{equation}
                 |d|< \sqrt{6N\log N} \quad \mbox{ and }\quad \left|t-\frac{2N}{3}\right|< \sqrt{2N\log N}. 
             \end{equation}

For positive integer  $m\leq n$  we have
\begin{eqnarray*}
                    P(D_N =2m, T_N =2n)&=& P(S_N =m+n, T_N =2n)  \\
                                       &=& \frac{N!}{3^{N}(n+m)! (n-m)! (N-2n)!} .
\end{eqnarray*}
\noindent
If  $d=2m$ and $t=2n$ satisfy $(2.1)$, then $(n+m)$, $(n-m)$ and  $(N-2n)$ are of exact order $N$ and we may apply Stirling's formula to see that
$$          P(D_N=2m,T_N =2n)$$ 
$$      = \frac{N^{N+\frac12}   \left(1+O\left(\frac 1N\right)\right)    }{ 2\pi  3^{N} (n+m)^{(n+m+\frac12)}(n-m)^{(n-m+\frac12)}(N-2n)^{(N-2n+\frac12)}}  $$
$$      = \frac{3^{\frac 32} \left(1+O\left(\frac 1N\right)\right) } { 2\pi N\left(\frac{ 3(n+m)}{N}\right)^{(n+m+\frac12)}\left( \frac{3(n-m)}{N}\right)^{(n-m+\frac12)} \left(\frac{3(N-2n)}{N}\right)^{(N-2n+\frac12)}} $$
$$      = \frac{3^{\frac32} \left(1+O\left( \frac 1N \right)\right)} { 2\pi N } \, \exp \Bigg\{ -\left(n+m+\frac12\right)\log\left( 1+\frac 3N \left(n+m-\frac N3\right) \right)  $$
$$ -\left(n-m+\frac12\right)  \log\left( 1+\frac 3N\left(n-m-\frac N3\right) \right)$$
$$ -\left(N-2n+\frac12 \right)\log\left(1+\frac 3N\left(\frac{2N}{3}-2n\right)\right)\Bigg\}.$$

\noindent
Next we expand the logarithms in the exponent. For the first order terms we obtain
$$-\frac 3N \bigg[\left(n+m+\frac 12\right)\left(n+m-\frac N3\right)+\left(n-m+\frac12\right)\left(n-m-\frac N3\right)+ $$
$$\left(N-2n+\frac 12\right)\left(\frac{2N}{3}-2n\right)\bigg]  $$ 
$$= -\frac 3N \left[\left(n+m-\frac N3\right )^2 + \left(n-m-\frac N3 \right)^2+\left( \frac{2N}{3}-2n\right)^2\right] $$
$$= -\frac{3}{N}\left( 6\,\tilde{n}^2+2m^2\right),$$
where $ \tilde{n}:= \left(n-\frac N3\right)$.\\

\noindent
The second order terms yield
$$ \frac 12 \left(\frac 3N\right)^2\bigg[ \left(n+m+\frac 12\right)\left(n+m-\frac N3\right)^2+\left(n-m+\frac 12\right)\left(n-m-\frac N3\right)^2$$
$$+\left(N-2n+\frac12\right)\left( \frac{2N}{3}-2n\right)^2 \bigg] $$
$$= \frac12 \left(\frac 3N \right)^2 \left[ -6\tilde{n}^3+(2N+3)\tilde{n}^2+6\tilde{n}m^2 +\left( \frac{2N}{3}+1 \right)m^2 \right]  $$
$$= \frac 3N \left( 3\tilde{n}^2+m^2 \right) + \frac{27}{N^{2}}\left(-\tilde{n}^3+\tilde{n}m^2 \right) +O\left( \frac{ \tilde{n}^2+m^2}{ N^{2}}\right).$$
\noindent
The third order terms contribute
$$-\frac 13\left(\frac 3N \right)^3 \bigg[\left(n+m+\frac12\right)\left(n+m-\frac N3\right)^3+\left(n-m+\frac12\right)\left(n-m-\frac N3\right)^3$$
$$+\left(N-2n+\frac12\right)\left(\frac{2N}{3}-2n\right)^3\bigg]$$
$$= \frac{18(\tilde{n}^3-\tilde{n}m^2)}{N^{2}}+O\left(\frac{\tilde{n}^4+m^4}{N^{3}}\right) $$ 
\noindent
As $d=2m$ and $t=2n$ satisfy $(2.1)$, the contribution of the remaining terms is dominated by that of the fourth order terms and equals
$$O\left(\frac{\tilde{n}^4+m^4}{N^{3}}\right). $$

Collecting the results of these computations we arrive at \smallskip

\noindent
\begin{equation} 
         \begin{split}
           & \qquad \qquad \qquad P(D_N=2m,T_N =2n) =  \frac{3^{\frac 32}}{  2\pi N}\,\\ 
                 &\times \exp\bigg\{ -\frac{3(3\tilde{n}^2+m^2)}{N}-\frac{9(\tilde{n}^3 
            -\tilde{n}m^2)}{N^2}  +O\left(\frac1N+ \frac{\tilde{n}^4+m^4}{N^3}\right) \bigg\} , 
\end{split} 
\end{equation}

\noindent
provided $m\leq n$ are integers between $1$ and $\frac12 N$ satisfying $m<\sqrt{2N\log N}$ and 
$|\tilde n|=\left|n-\frac N3\right|<\sqrt{N\log N}$. However, we shall also use $(2.2)$ if these inequalities do not hold, since in that case both left- and right-hand members of $(2.2)$ are negligible for our purposes.

By Taylor expansion of the integrand about $x=m$, we find that for integer $0<m\leq n$ with $m<\sqrt{2N\log N} $ and $|\tilde n| = \left|n-\frac N3 \right|<\sqrt{N\log N}$,
\begin{eqnarray*}
           \int_{[m-\frac12,m+\frac12) } \exp\bigg\{ -\frac{3}{N}x^2 +\frac{9}{N^2}\tilde{n}x^2  
                                           +O\left( \frac1N+\frac{1}{N^3}\left(\tilde{n}^4+x^4\right)\right) \bigg\}dx \\[3mm]
        = \exp\bigg\{-\frac{3}{N}m^2 + \frac{9}{N^2}\tilde{n}m^2+  
        O\left(\frac1N+ \frac{1}{N^3}\left(\tilde{n}^4+m^4\right)\right)\bigg\}.      
\end{eqnarray*}
\noindent
It follows that for integers $0<m\leq n$ with $m<\sqrt{2N\log N} $ and $|\tilde n| = \left|n-\frac N3 \right|<\sqrt{N\log N}$,
\begin{eqnarray*}
&&  P(2\leq D_N\leq 2m,T_N =2n) \\ 
&&=\frac{ 3^{\frac32}}{ 2\pi N}\int_{[\frac12,m+\frac12)} e^{-\left\{\frac{3}{N}x^2 -\frac{9}{N^2}\tilde{n}x^2 +O\left(\frac{1}{N^3}x^4\right)\right\}}dx \,
e^{-\left\{ \frac{9}{N}\tilde{n}^2 +\frac{9}{N^2}\tilde{n}^3  +O\left(\frac1N+ \frac{1}{N^3}\tilde{n}^4 \right)\right\} }\\ 
&& = \frac{3^{\frac32}}{ 2\pi N}\int_{[\frac 12,m+\frac12)} e^{-\left\{\frac{3}{N}x^2 -\frac{9}{N^2}\tilde{n}x^2 \right\}}dx \, 
 e^{-\left\{ \frac{9}{N}\tilde{n}^2 + \frac{9}{N^2}\tilde{n}^3  + O\left(\frac 1N+\frac{1}{N^3}\tilde{n}^4\right)\right\}}, 
\end{eqnarray*}
\noindent
and again we may use this for all integers $m$ and $n$ with $0<m\leq n\leq \frac12 N$ with impunity.

For real $r>\frac12 $ we write $r=m+\theta$ where $m=\lfloor r\rfloor$ and $\theta= \Frc (r)=r-\lfloor r \rfloor \in[0,1)$ denote the integer and fractional parts of r respectively. Then for $r<\sqrt{2N\log N}$ and $|\tilde n| = \left|n-\frac N3\right|<\sqrt{N\log N}$,
\begin{eqnarray*}
 P(2\leq D_N\leq 2r,T_N =2n) = P(2\leq D_N \leq 2m,T_N =2n)=  \qquad \\[3mm] 
 \frac{3^{\frac32}}{ 2\pi N}e^{-\left\{\frac{9\tilde{n}^2}{N}+ \frac{9\tilde{n}^3}{N^2}+ O\left(\frac1N+\frac{\tilde{n}^4}{N^3}\right)\right\}}  \bigg[ \int\limits_{[\frac12,r)}e^{\frac{-3x^2}{N}+\frac{9\tilde{n}x^2}{N^2}}dx  
+   \int\limits_{[r,m+\frac12)} e^{\frac{-3x^2}{N}+\frac{9\tilde{n}x^2}{N^2}} dx \bigg]  .
\end{eqnarray*}
\noindent
Evaluating the second integral by expanding the integrand about the point $x=m+\frac12 $, we arrive at
\begin{eqnarray*}
&& P(2\leq D_N\leq 2r,T_N =2n)  = 
\frac{3^{\frac32}}{ 2\pi N}e^{-\left\{ \frac{9}{N}\tilde{n}^2 +\frac{9}{N^2}\tilde{n}^3+
 O\left(\frac1N+\frac{1}{N^3}\tilde{n}^4  \right)\right\}} \\ && \quad  \quad\times 
 \left[ \int_{[\frac12,r)} e^{-\frac{3}{N}x^2 +\frac{9}{N^2}\tilde{n}x^2}dx -
 e^{-\frac{3}{N}r^2 +\frac{9}{N^2}\tilde{n}r^2 } \left(\Frc (r)-\frac12+O\left(\frac rN\right)\right) \right]. 
\end{eqnarray*}
\noindent
Again we may use this for all $r>0$ and integer $n\leq\frac12 N$.

Choose $w>0$ and $r=w\sqrt{\frac n2}$. We have
\begin{eqnarray*}
&& P(0<W_N\leq w, T_N \mbox{ is even }) = \sum_{1\leq n\leq \frac N2 } P(2\leq D_N\leq 2r, T_N=2n)\\ 
&&= \sum_{1\leq n\leq \frac N2} \frac{3^{\frac32}}{2\pi N}e^{-\left\{\frac{9}{N}\tilde{n}^2  +\frac{9}{N^2}\tilde{n}^3
+O\left(\frac1N+\frac{1}{N^3}\tilde{n}^4\right)\right\}}\\
 && \qquad \times  
\left[ \int_{[\frac12,r)}e^{-\frac{3}{N}x^2 +\frac{9}{N^2}\tilde{n}x^2 }dx-
e^{-\frac{3}{N}r^2 +\frac{9}{N^2}\tilde{n}r^2 }  \left( \Frc (r)-\frac12 +O\left(\frac rN \right)\right)         \right].  
\end{eqnarray*}

\noindent
As $|\tilde n| = |n- \frac N3|\leq \frac N3$, the expression between square brackets is of order $\sqrt{N}$. Next, comparison with the normal $N(\frac N3,\frac {N}{18})$ distribution with mean 
$\frac N3$ and variance $\frac {N}{18}$ shows that  

\begin{equation*}
\frac{3}{\sqrt{\pi N}} \sum_{1\leq n\leq \frac N2 } e^{-\left\{\frac{9}{N}\tilde{n}^2 +\frac{9}{N^2}\tilde{n}^3+ O\left(\frac1N+\frac{1}{N^3}\tilde{n}^4 \right)\right\}}      
   = \frac{3}{\sqrt{\pi N}} \sum_{1\leq n\leq \frac N2} e^{-\frac{9}{N}\tilde{n}^2 } + O\left(\frac 1N\right),  
\end{equation*}
\noindent
and hence 

\noindent
        \begin{equation}
                \begin{split}
                 P(0<W_N \leq w, T_N \mbox{ is even}) 
                 = O\left( \frac 1N\right) + \frac12 \sum_{1 \leq n \leq \frac N2 }\frac{3}{\sqrt{\pi N}}e^{-\frac{9}{N}\tilde{n}^2 }\sqrt{\frac{3}{\pi N}} 
\\ \times \left[ \int\limits_{[\frac12 ,r) }e^{-\frac{3}{N}x^2 +\frac{9}{N^2}\tilde{n}x^2  }dx- 
                e^{-\frac{3}{N}r^2 +  \frac{9}{N^2}\tilde{n}r^2  }   \left(\Frc (r)-\frac12 +O\left(\frac rN\right)\right)\right] \\ 
             \end{split}
\end{equation}

\noindent
Since $\sqrt{ \frac{3}{\pi N} }e^{ -\frac{3}{N}x^2}$ is the density of the $N\left(0,\frac N6\right)$ distribution, we see that
$$ \sqrt{ \frac{3}{\pi N} } \int_{[\frac12,r)} e^{-\frac{3}{N}x^2+\frac{9}{N^2}\tilde{n}x^2 }dx =  
\sqrt{ \frac{3}{\pi N} } \int_{[\frac12,r)} e^{-\frac{3}{N}x^2}dx + \frac{3 \tilde{n} g(r)}{2N} + O\left(\frac 1N\right),  $$
\noindent
where
$$ g(r)= \frac{1}{\sqrt{2\pi} } \int_B x^2 e^{-\frac12 x^2}dx,  $$\smallskip
\noindent
with $B=\left( \frac12\sqrt{\frac 6N},r\sqrt{\frac 6N} \right)$. Now 
$$r\sqrt{\frac 6N}=w\sqrt{\frac{3n}{N}} = w+\frac{3w\tilde{n}}{2N}+O\left(w\left(\frac{\tilde{n}}{N}\right)^2\right), $$
and splitting the integral in one over $\left(\frac12 \sqrt{\frac 6N},w\right)$ and one over $\left(w,r\sqrt{\frac 6N} \right)$ and expanding the latter around the point $w$, we obtain
$$g(r)= h(w,N)+O\left(\frac{\tilde{n}}{N} \right),$$
for a bounded function $h$. It follows that

$$ \sqrt{ \frac{3}{\pi N} }  \int_{[\frac12 ,r)} e^{-\frac{3}{N}x^2 +\frac{9}{N^2}\tilde{n}x^2}dx =   $$

$$ =  \sqrt{ \frac{3}{\pi N} } \int_{[\frac12,r)} e^{-\frac{3}{N}x^2}dx + \frac{3\tilde{n}\, h(w,N)}{2N} + O\left( \frac 1N+ \left(\frac{\tilde{n}}{ N}\right)^2 \right).  $$
Substituting this in $(2.3)$ and comparing the distribution of n once more with $N\left( \frac N3,\frac{N}{18} \right)$, we see that the last two terms above contribute only $O(N^{-1})$. Finally the term containing $N^{-\frac12}e^{-\frac{3}{N}r^2}$. $O\left(\frac rN\right)$ is clearly $O(N^{-1})$. Hence we have reduced $(2.3)$ to 
\begin{equation}
 \begin{split}
         P(0<W_N\leq w, T_N \mbox{ is even}) = 
          \frac12 \sum_{1\leq n\leq \frac N2} \frac{3}{\sqrt{\pi N}} e^{-\frac{9}{N}\tilde{n}^2 } \sqrt{\frac{3}{\pi N}} \\[3mm]
          \times [ \int_{[\frac12 ,r)} e^{-\frac{3}{N}x^2}dx -e^{-\frac{3}{N}r^2+\frac{9}{N^2}\tilde{n}r^2 } 
        \left( \Frc (r)-\frac12 \right) +O\left(\frac 1N \right).  
\end{split}
\end{equation}
\noindent
Now
$$ \frac{1 }{\sqrt{N}}e^{-\frac{3}{N}r^2 +\frac{9}{N^2}\tilde{n}r^2 } = \frac{1}{\sqrt{N}}e^{-\frac12 w^2n \left(\frac3N-\frac{9}{N^2}\tilde{n} \right)} =  \frac{1}{\sqrt{N}} e^{ -\frac12 w^2\left( 1- \frac{9}{N^2}\tilde{n}^2\right) }$$
$$ =\frac{1}{\sqrt{N}}e^{-\frac12w^2}\left(1+O\left( \frac{9}{N^2}\tilde{n}^2\right)\right) $$
and the remainder term will give rise to another $O(N^{-1})$ term in $(2.4)$. We obtain
$$  P(0<W_N\leq w, T_N \mbox{ is even}) =O\left(\frac 1N\right) +$$
$$ \frac12 \sum_{1\leq n\leq \frac N2} \frac{3}{\sqrt{\pi N}} e^{-\frac{9}{N}\tilde{n}^2 } \sqrt{\frac{3}{\pi N}}   \left[ \int_{[\frac12,r)} e^{-\frac{3}{N}x^2 }dx- e^{-\frac12 w^2} 
\left\{\Frc\left(w\sqrt{\frac{n}{2}}\right) - \frac 12\right\} \right]. $$
\medskip

As $r=w\sqrt{n/2}$ ,we have
$$ \sqrt{\frac{3}{\pi N}}  \int_{[\frac12,r)} e^{-\frac{3}{N}x^2}dx = \Phi\left(w\sqrt{\frac{3n}{N}} \right) -\Phi\left( \sqrt{ \frac {3}{2N} } \right) $$
$$= \Phi(w)-\frac12 - \sqrt{ \frac{3}{4\pi N} }+ \frac{3\tilde{n}w\varphi (w)}{2N}+ O\left( \left(\frac{\tilde{n}}{N}\right)^2 \right),  $$
where $\Phi$  and $\varphi$ denote the standard normal distribution function and its density.  Obviously, the linear term containing $\tilde{n}$ as well as the remainder term contribute $O(N^{-1})$ to $(2.4)$. Since
$$ \sum_{1\leq n\leq \frac N2} \frac{3}{\sqrt{\pi N}} e^{-\frac{9}{N}\tilde{n}^2 }= 1+O\left(\frac1N\right) ,  $$
we find 
\begin{equation}
         \begin{split}
          &       P(0<W_N \leq w, T_N \mbox{ is even}) =  
\frac12 \left(\Phi (w)-\frac12 -\sqrt{ \frac{3}{4\pi N}}\right) \\
 &\quad -\sqrt{\frac{3}{2N}}\varphi(w)\sum_{1\leq n\leq \frac N2} \frac{3}{\sqrt{\pi N}} e^{-\frac{9}{N}\tilde{n}^2 }      \left\{\Frc\left(w\sqrt{\frac{n}{2}}\right) - \frac 12\right\} +O\left(\frac 1N \right) 
              \end{split}
\end{equation}

  An almost identical computation produces an asymptotic expression for 
the case where $T_N$ is odd. We have 
\begin{equation}
        \begin{split}
                &P(0<W_N \leq w, T_N \mbox{ is odd} ) = \frac 12 \left(\Phi(w)-\frac12 \right)  - \sqrt{\frac{3}{2N}} \varphi(w)\  \times \\
                &\sum_{0\leq n\leq \frac{N-1}{2}} \frac{3}{\sqrt{\pi N}} e^{-\frac{9}{N}(n+\frac12 -\frac N3)^2 } \left\{\Frc\left(\frac{1}{2}w\sqrt{2n+1} + \frac{1}{2}\right) - \frac 12\right\} +O\left(\frac1N \right). 
        \end{split}
\end{equation}
A few minor modifications of $(2.5)$ and $(2.6)$ are in order. In both formulas we may replace the sum by a sum over $0 \leq n \leq N$ since the terms added are exponentially small in $N$. Next, in $(2.6)$,  we may change the factor $\left(n+\frac12 -\frac N3\right)^2$ in the exponent to $\tilde n^2 = (n- \frac N3)^2$ because the additional factor $e^{O\left( (n-\frac N3)/N \right)}=1+O\left(\tilde n/N\right)$ only produces an additional remainder term of order $N^{-1}$. Also in $(2.6)$, we may replace $\Frc \left( \frac 12 w \sqrt{2n+1}+\frac12\right)$ by $\Frc\left(\frac12 w \sqrt{2n}+\frac12 \right) =\Frc\left( w\sqrt{\frac n2}+\frac12\right)$. This is because $(\frac 12 w\sqrt{2n+1} + \frac 12 ) - (\frac 12 w\sqrt{2n} + \frac 12)$ is of exact order $\frac{w}{\sqrt n}$ and hence $\frac{w}{\sqrt N}$, so $\Frc\left( \frac12 w \sqrt{2n+1} +\frac12 \right) - \Frc\left( \frac12 w \sqrt{2n}+\frac12\right)$ is roughly equal to $-1$ for only one out of every $\frac{c\sqrt{N}}{w}$ values of $n$, i.e.\ for $O\left(w\sqrt{N}\right)$ values of $n$ at distances of exact order $\sqrt{N}/w$. For all other $n$, this difference is $O(N^{-\frac12})$ and this implies that the substitution only adds another $O(N^{-1})$ remainder term. Hence
\begin{equation}
        \begin{split}
                P\left(0 < W_N\leq w, T_N \mbox{ is odd} \right) = \frac 12 \left(\Phi(w)- \frac12\right) - \sqrt{ \frac{3}{2N}} \varphi(w) \\
\times  \sum_{0\leq n\leq N} \frac{3}{\sqrt{\pi N}} e^{-\frac 9N \tilde n^2}\left[\Frc\left(w\sqrt{\frac n2}  +\frac 12\right)-\frac12 \right]+O\left(\frac{1}{N}\right).
        \end{split}
\end{equation}
If $0\leq \Frc\left( w\sqrt{\frac n2} \right)<\frac12$, then $\Frc\left( w\sqrt{\frac n2}\right)+\Frc\left(w\sqrt{\frac n2}+\frac12\right) = 2\Frc\left( w\sqrt{\frac n2}\right)+\frac12 = \Frc\left(2w\sqrt{\frac n2}\right)+ \frac12 $. On the other hand, if $\frac12 \leq \Frc\left( w\sqrt{\frac n2}\right)<1$, then $\Frc\left( w\sqrt{\frac n2}\right) +\Frc\left( w\sqrt{\frac n2}+\frac12 \right) = 2\Frc\left( w\sqrt{\frac n2}\right)-\frac12 = \Frc\left( 2w\sqrt{\frac n2}\right)+ \frac12 $ . Hence $\Frc\left( w\sqrt{\frac n2}\right) +\Frc\left( w\sqrt{\frac n2}+\frac12\right) = \Frc \left( 2w\sqrt{\frac n2}\right)+ \frac12 = \Frc\left( w\sqrt{2n} \right)+ \frac12$  in both cases. Hence we may combine $(2.5)$ and $(2.7)$ to obtain for $w>0$,
\begin{equation}
        \begin{split}
                P\left( 0<W_N\leq w \right) = \Phi(w)-\frac12 -\frac12 \sqrt{ \frac{3}{4\pi N}} - \sqrt{\frac{3}{2N}} \varphi(w) \\
\times                \sum_{ 0\leq n\leq N}  \frac{3}{\sqrt{\pi N}} e^{-\frac 9N \left(n-\frac{N}{3}\right)^2 } \left[ \Frc\left(w \sqrt{2n}\right)-\frac12 \right]+O\left( \frac{1}{N} \right) .
        \end{split}
\end{equation}
 
Because the distribution of $W_N$ is symmetric about the origin and $P(W_N=0)= \sqrt{\frac{3}{4\pi N}}+ O(N^{-1})$, this determines the expansion for the distribution function $F_N$ of $W_N$. The expansion is uniform in $w>0$. Since it is identical to $(1.6)$, this proves the third statement of the theorem.

It remains to prove that any closed interval of length $O(N^{-1})$ that does not contain the origin has probability $O(N^{-1})$  under $\Psi_N$ and hence $F_N$. Clearly, this will imply the second statement of the theorem. Obviously, the only term in $(2.8)$ that we need to consider is 
\begin{equation*}
R(w) = \frac{\Lambda_N(w)}{\sqrt{N}}
\end{equation*}
\begin{equation}
= -\sqrt{\frac{3}{2N}}  \varphi (w) \sum_{0\leq n\leq N} \frac{3}{\sqrt{\pi N} }e^{-\frac 9N\left(n - \frac N3\right)^2}  
        \left(\Frc \left(w\sqrt{2n} \right)-\frac12 \right) ,
\end{equation}
as the remainder of the expansion obviously has bounded derivative. 

 We begin by noting that if for a given $w>0$, $w\sqrt{ 2n}$ is an integer for some $1\leq n\leq  N$, then $\Frc \left(w\sqrt{ 2n}\right)$ and hence $R$ has a jump discontinuity at this value of $w$. In the range where
$|n- \frac N3|= x\sqrt{N}$ for $|x|\leq y$, there can be a most $wy$ such integer values of $n$. To see  this, simply note that if $w\sqrt{ 2n}=k$ and $w\sqrt{2n'}=k+1$ , then $|n'-n| \geq \frac{ 2\sqrt{N}}{w}$ , so there can be only 
$\frac{2y\sqrt{N}}{ \frac{2\sqrt{N}}{w} }=wy$ values of $n$ in the required interval. Such a value of n contributes an amount $O\left(N^{-1}\varphi(w)e^{-9x^2}\right)$ to the jump discontinuity at $w$, and hence $R(w)-R(w-0)= O(N^{-1})$ at such a point $w$. Incidentally, this proves the second part of Theorem 1.1.

 Choose $\epsilon>0$ and consider two such jump points $w\not=w'$ in $[\epsilon , \infty)$ with $w\sqrt{  2n}=k$ and $w'\sqrt{2n'}=k'$ for integers $k, k', n$ and $n'$ with $(n-\frac N3)=x\sqrt{N},\ \left(n'-\frac N3\right)=x'\sqrt{N} $ and $|x|\vee | x'|\leq y$. Suppose that $(w'-w)=O(N^{-1})$ and hence $\frac{w'-w}{w}=O(N^{-1})$ since $w\geq \epsilon$. For given $w$, $n$ and $k$, we ask how many integer values of $n'$ satisfy these conditions.

First we note that, for some positive $c$ there are only at most $cw(y+1)$ possible choices for $k'$ since $\sqrt{2n} = \sqrt{2\frac N3 +2x\sqrt{N}} = \sqrt{\frac{2N}{3}} +\sqrt{\frac32}x+O\left( \frac{y^2}{\sqrt{N}} \right), \sqrt{2n'} =\sqrt{\frac23 N} +\sqrt{\frac 32}x'+O\left(\frac{y^2}{\sqrt{N}}\right)$ and hence $|k'-k|\leq 2wy + O\left( w\frac{y^2}{\sqrt N}+|w'-w|\sqrt{N}\right)\leq \left(\frac c2 \right)w(y+1)$. For each choice of $k'$, the corresponding $n'$ satisfies $n' =\frac12 \left(\frac{k'}{w'}\right)^2$ for some admissible $w'$, and since $w,w'\geq \epsilon$ and $(w'-w)=O(N^{-1})$, this leaves a range of order $O\left( \left(\frac{k'}{w'}\right)^2 N^{-1} \right)=O(1)$ for $n'$. Hence, for some $C>0$, there are at most $Cw(y+1)$ possible values of $n'$ for which there exists an integer $k'$ with $(w'-w) = O(N^{-1})$. By the same argument as above, the total contribution of discontinuities to$|R(w')-R(w) |$ is $O(N^{-1})$ as long as $|w-w'|= O(N^{-1})$. As any closed interval of length $O(N^{-1})$ that does not contain the origin is bounded away from $0$, this holds for the sum of the discontinuities in such an interval.

 At all other points $w>0$, $R$ is differentiable and the derivative of $\Frc \left(w\sqrt{ 2n}\right)$ equals $\sqrt{2 n}$. Hence the derivative of $R$ is $O(1)$ and its differentiable part contributes at most $O(N^{-1})$ to the probability of any interval of length $O(N^{-1})$. This completes the proof of the Theorem 1.1.

\section{Evaluation of the oscillatory term}

Let $W$ denotes a r.v. with non negative c.f. $\psi(t)\ge 0$ of
support contained in $[-1,1]$ and exponential decay of density 
of type $\exp\{-|x|^{2/3}\},\,x\to\infty$, 
\citep[see e.g.] [p. 85] {Bhattacharya-Rao:1986/2}.
Introduce  r.v. $w_N := w +N^{-1/2}(\log N)^{-1}\4 W, \, w >0$ and 
let $c>0$ denote  an positive absolute constant.   
Then we  may bound   the normal approximation
error in $(1.6)$ using  similar arguments as in the proof of  
the well-known smoothing inequality, (see Lemma 12.1 of  Bhattacharya and Rao),
 obtaining, for $w\ge 1$,
\begin{equation}
\begin{split}
\,N^{-1/2}\bigl|\Ee \Lambda_N(w_N)\bigr| \le & 
\, \bigl|\Ee \bigl(F_N(w_N) -\Phi(w_N)\bigr)\bigr| +
cN^{-1}\\ 
& \,\le \sup_{x\in[w-1/2,w+1/2]}\bigl|F_N(x) -\Phi(x)\bigr|+cN^{-1}, 
\end{split}
\end{equation} 
where
$$\Lambda_N(w):=  \, -\varphi(w)\4 
\sum_{1\le n\le N}\frac {3^{3/2}}{\sqrt{2 \pi N}}
\exp\{-\frac 9 N(n - \frac N 3)^2\}(\Frc(w \sqrt {2n})- 1/2 ).$$
We start with the following Fourier series expansion
\beqa 
\tau(x):= frac(x)-1/2 = -\sum_{k=1}^{\infty} 2\4  \frac{ sin(2 \pi\4 k\4 x)}
{2\4 k\4 \pi},  
\eeqa
which holds for all nonintegral $x$. 

Note that by the properties of $W$ (i.e. the vanishing of Fourier coefficients)
$$ 
\Ee \tau(w_N\sqrt{2n})=  -\sum_{k=1}^{M_n}  \Ee \frac{ sin(
  2\pi\4 k\4 \sqrt{2n}(w+N^{-1/2}(\log N)^{-1}\4 W))} { k\4 \pi}, 
$$
where   $M_n:=[\sqrt N\log N/(2\pi\sqrt {2n})]+1$, i.e. $M_n = O(\log N)$ for 
$|n-N/3|<\sqrt{N\log N}$.

Rewriting $\Lambda_N(w)$ in $(3.1)$ in the form
\begin{equation}
\Lambda_N(w) := \, -\frac {3^{3/2}} {(2 \4 \pi\4 N)^{1/2}} \varphi(w) 
\sum_{n=1}^{N} exp\{-9\4 (\tilde{n}^2/N\}\4 
\tau\bigl(w\4(2 \4 n)^{1/2}\bigr), 
\end{equation}
where $\tilde{n}:=n-N/3$,
we get 
\begin{equation}
\begin{split}
\Ee \Lambda_N(w_N) =& \, \sqrt{\frac 3 2}\pi^{-1} 
\sum_{k=1}^{M} \frac 1 k 
\lambda_{N,k}+O(N^{-3}), \,\, 
\text{where}\,\, M:=[\log N]\,\, \text{and}\\
\lambda_{N,k}  :=& \, \frac 3{\sqrt{\pi \4 N}} \Ee \varphi(w_N)\sum_{n=1}^{N} 
exp\{-9\4 
\tilde{n}^2/N\}\4  sin(2\pi\4 k\4 w_N\4  \sqrt{2n}). 
\end{split}
\end{equation}
In the arguments of the $\sin$ function we
 use a Taylor expansion, for $|n-N/3|<\sqrt{N\log N}$, 
$$
\sqrt{n} = \sqrt{N/3} + \sqrt 3 \tilde{n}/(2\sqrt N)+
O\bigl(\tilde{n}^2/N^{3/2} \bigr).
$$
Thus, for $|\tilde{n}|<\sqrt{N\log N}$,  
\begin{equation}
 sin(2\pi\4 k\4 w_N\4 \sqrt{2n})=  sin\bigl(d_0 + \4 \pi\4 d_1 \tilde{n} \bigr) 
+  O\bigl(k\4 w_N N^{-3/2}\4 \tilde{n}^2\bigr), 
\end{equation}
where $ d_0 :=2\pi\4k \4 w_N\4 (\frac 2 3)^{1/2} \sqrt{N}$, 
$d_1 := k\4 \4 w_N\4
(\frac 3 2)^{1/2}\4/ \sqrt N$. %
Hence  we may write
\begin{equation}
\begin{split}
 \lambda_{N,k} =& \frac 3 {\sqrt{\pi \4 N}} \Ee \varphi(w_N) 
\sum_{n\in \Ze} exp\{-9\4 
\tilde{n}^2/N\}\4
\sin\bigl(d_0 + \4 2\pi \4d_1 \4\tilde{n}\bigr)\\
&\q  + O(k  N^{-1/2} \4 \log N ).  
\end{split}
\end{equation}
We shall now evaluate the theta sum on the left hand side using
Poisson's formula, 
\citep[see e.g.][p. 189]{Mumford:1983}.
\begin{equation}
\sum_{m \in \Ze} \exp\{-z\4 m^2 +i2\pi\4 m \4 b\} = \pi^{1/2}z^{-1/2}
\sum_{l \in \Ze}
\exp\{- \pi^2 \4 z^{-1}(l- b)^2\}, 
\end{equation}
where $b \in \mathbb R$, $\Re z >0$ and $z^{1/2}$ denotes the branch with 
positive real part.
Writing $sin(x)=(\exp[i\4 x] - \exp[-i\4 x])/2$ in (3.6) 
and assuming for simplicity $N/3 \in \Ze$ we may replace  summation over n by
summation over $ m:=\tilde{n}= n - N/3 \in \Ze$ in (3.6). 
Applying now (3.7) we have  to bound the imaginary part of 
 expectations of  theta functions of type 
\beqa
I_k:=\frac 3 {\sqrt{\pi\4 N}} \Ee \varphi(w_N)\exp\{\4 i\4 d_0 \}
\sum_{m\in \Ze} \exp\{-
9\4 m^2 \4 N^{-1} + i\4 2\pi\4 d_1\4m\}. 
\eeqa
We obtain for $ k \le M =[\log N]$ that $|d_1| \le 2
 \4N^{-1/2}(\log N) |w_N| \le 4 \4N^{-1/2}(\log N)^2 $ with probability
 $1- O(N^{-3/2})$ by the assumption $w \le \log N$.
   Hence the dominant term
in (3.9) below  is the term with $l=0$ and we obtain with 
$c_{N,k}:=  \exp\{ \4 2\pi\4 i\4 \4k \4 w_N\4 (\frac 2 3)^{1/2} \sqrt{N}\}$ 
\begin{equation}
\begin{split}
I_k   =& \,  \Ee c_{N,k}\varphi(w_N)\4 \sum_{l \in \Ze} 
\exp\{ - N \4( l- d_1)^2\4 \pi^2/9\}
\\
   = & \,  \Ee  c_{N,k}\4 \varphi(w_N)\exp\{-  N\4 d_1^2\4 \pi^2/9\} + O\bigl(N^{-3/2}\bigr)
\\ 
 =& \,\4 f_{N,k}\4 \varphi(w) \exp\{- \pi^2 k^2 \4  w^2/6  
+ \4 i\4 \42\pi k \4 w\4 (\frac 2 3)^{1/2} \sqrt{N}\}
+  O\bigl(N^{-1/2}(\log N)^4\bigr), 
\end{split}
\end{equation}
where $f_{N,k}:= \psi\bigl(2\pi\4 (\frac 2 3)^{1/2} \4 \frac k {\log
  N}\bigr)= 1+ O\bigl((k/\log N)^2\bigr)$. 
Using the equation (3.9) in (3.4)  we get
\begin{equation}
\begin{split}
\Ee \Lambda_N(w_N) =& \, \sqrt\frac 3 2 
\varphi(w)\Im \sum_{k=1}^M \frac {f_{N,k}} {k\pi}
\exp\bigl\{- \frac{\pi^2}6\4 k^2\4  w^2\4+ 
2\pi\4 i\4 \4k \4 w\4 \sqrt{\frac{2N}3}\bigr\}\\
& \, \q + O\bigl(N^{-1/2}(\log N)^5\bigr). 
\end{split}
\end{equation}
Hence, there exists a constant $c_0(w) >0$ such that
\begin{equation}
|\Ee \Lambda_N(w_N)| > c_0(w) >0,
 \end{equation} 
provided that  $4\4 w\4 \sqrt{\frac{2 N} 3} $ is an odd integer, which
proves the assertion $(1.10)$. 

\renewcommand\bibsection{\section*{REFERENCES}}
\bibliographystyle{ims}
\bibliography{Bickel-datenbank1}

\begin{thebibliography}{11}
\expandafter\ifx\csname natexlab\endcsname\relax\def\natexlab#1{#1}\fi
\expandafter\ifx\csname url\endcsname\relax
  \def\url#1{\texttt{#1}}\fi
\expandafter\ifx\csname urlprefix\endcsname\relax\def\urlprefix{URL }\fi
\providecommand{\eprint}[2][]{\url{#2}}

\bibitem[{Albers et~al.(1976)Albers, Bickel and van
  Zwet}]{Albers-Bickel-vanZwet:1976}
\textsc{Albers, W.}, \textsc{Bickel, P.~J.} and \textsc{van Zwet, W.~R.}
  (1976).
\newblock Asymptotic expansions for the power of distribution free tests in the
  one-sample problem.
\newblock \textit{Ann. Statist.}, \textbf{4} 108--156.

\bibitem[{Bentkus and G\"{o}tze(1996)}]{Bentkus-Goetze:1996}
\textsc{Bentkus, V.} and \textsc{G\"{o}tze, F.} (1996).
\newblock The {B}erry-{E}sseen bound for {S}tudent's statistic.
\newblock \textit{Ann. Probab.}, \textbf{24} 491--503.

\bibitem[{Bentkus et~al.(1997)Bentkus, G{\"o}tze and van
  Zwet}]{bentkus-goetze-vanzwet:1997}
\textsc{Bentkus, V.}, \textsc{G{\"o}tze, F.} and \textsc{van Zwet, W.~R.}
  (1997).
\newblock An {E}dgeworth expansion for symmetric statistics.
\newblock \textit{Ann. Statist.}, \textbf{25} 851--896.

\bibitem[{Berry(1941)}]{Berry:1941}
\textsc{Berry, A.~C.} (1941).
\newblock The accuracy of the {G}aussian approximation to the sum of
  independent variates.
\newblock \textit{Trans. Amer. Math. Soc.}, \textbf{49} 122--136.

\bibitem[{Bhattacharya and Ranga~Rao(1986)}]{Bhattacharya-Rao:1986/2}
\textsc{Bhattacharya, R.} and \textsc{Ranga~Rao, R.} (1986).
\newblock \textit{Normal Approximation and Asymptotic Expansions}.
\newblock Wiley, New York.

\bibitem[{Bickel et~al.(1986)Bickel, G\"otze and {van
  Zwet}}]{Bickel-Goetze-vanZwet:1986}
\textsc{Bickel, P.~J.}, \textsc{G\"otze, F.} and \textsc{{van Zwet}, W.~R.}
  (1986).
\newblock The {E}dgeworth expansion for {$U$}-statistics of degree two.
\newblock \textit{Ann. Stat.}, \textbf{14} 1463--1484.

\bibitem[{Brown et~al.(2002)Brown, Cai and DasGupta}]{Brown-Cai-DasGupta:2002}
\textsc{Brown, L.~D.}, \textsc{Cai, T.~T.} and \textsc{DasGupta, A.} (2002).
\newblock Confidence intervals for a binomial proportion and asymptotic
  expansions.
\newblock \textit{Ann. Statist.}, \textbf{30} 160--201.

\bibitem[{Esseen(1942)}]{Esseen:1942}
\textsc{Esseen, C.-G.} (1942).
\newblock On the {L}iapounoff limit of error in the theory of probability.
\newblock \textit{Ark. Mat. Astr. Fys.}, \textbf{28A} 19.

\bibitem[{Feller(1965)}]{Feller:1965/2}
\textsc{Feller, W.} (1965).
\newblock \textit{An introduction to probability theory and its applications},
  vol.~2.
\newblock Wiley, New York.

\bibitem[{Mumford(1983)}]{Mumford:1983}
\textsc{Mumford, D.} (1983).
\newblock \textit{Tata lectures on theta. {I}}, vol.~28 of \textit{Progress in
  Mathematics}.
\newblock Birkh\"auser Boston Inc., Boston, MA.
\newblock With the assistance of C. Musili, M. Nori, E. Previato and M.
  Stillman.

\bibitem[{van Zwet(1984)}]{vanZwet:1984}
\textsc{van Zwet, W.~R.} (1984).
\newblock A {B}erry-{E}sseen bound for symmetric statistics.
\newblock \textit{Z. Wahrsch. Verw. Gebiete}, \textbf{66} 425--440.

\end{thebibliography}

\end{document}